\documentclass[12pt, a4paper]{article}

\usepackage[a4paper, margin=1in]{geometry}
\usepackage[utf8]{inputenc}
\usepackage[T1]{fontenc}
\usepackage{lmodern}
\usepackage{amsmath,amssymb,amsfonts,amsthm}
\usepackage{mathtools}
\usepackage{upgreek}
\numberwithin{equation}{section}

\usepackage{graphicx}
\usepackage{caption}
\usepackage{subcaption}
\usepackage{booktabs}
\usepackage{multirow}
\usepackage{enumitem}

\usepackage[round,authoryear]{natbib}
\usepackage[colorlinks=true, linkcolor=blue, citecolor=blue, urlcolor=blue]{hyperref}

\theoremstyle{plain}
\newtheorem{theorem}{Theorem}[section]
\newtheorem{proposition}[theorem]{Proposition}

\newtheorem{corollary}[theorem]{Corollary}
\theoremstyle{definition}
\newtheorem{definition}[theorem]{Definition}

\theoremstyle{remark}
\newtheorem{remark}{Remark}

\begin{document}

\begin{center}
{\Large\bf A New Perspective on Cumulative Residual R\'enyi Entropy} \\
\vspace{0.3in}
{\bf Shubham Babasaheb Alte$^{a}$, Swaroop Georgy Zachariah$^{b}$ and Mohd. Arshad$^{a,}$\footnote{Corresponding author. E-mail: arshad@iiti.ac.in (M. Arshad).}} \\
\vspace{0.1in}
{\small $^a$ Department of Mathematics, Indian Institute of Technology Indore, India \\
$^b$ Department of Mathematical Sciences, Rajiv Gandhi Institute of Petroleum Technology, Bengaluru Campus, Karnataka, 562165, India}
\end{center}

\noindent\rule{\textwidth}{0.4pt}
\begin{abstract}
This paper introduces a new cumulative residual R\'enyi entropy (CRRE), which addresses limitations of the existing cumulative residual R\'enyi entropy. The conventional measure suffers from two primary drawbacks: its failure to converge to the cumulative residual entropy (CRE) as the order parameter approaches unity, and its lack of guaranteed positivity. Our proposed CRRE resolves these issues, offering a theoretically consistent and always nonnegative measure of uncertainty in residual lifetimes. We explore its mathematical properties, including its relationship with existing information measures. Furthermore, the concept is extended to its dynamic version, and its time-dependent characteristics are investigated. A non-parametric estimator of the proposed CRRE is developed, its almost sure convergence is established, and is validated through Monte Carlo simulation. The effectiveness of CRRE as a tool for quantifying uncertainty is confirmed through its application in chaotic theory and financial data.
\end{abstract}
\vskip 2mm

\noindent \emph{Keywords}: R\'enyi entropy; Tsallis entropy; cumulative residual entropy; stochastic ordering; empirical estimation.\\
\noindent \emph{Mathematics Subject Classification (2020)}: 94A17, 62B10, 62G05, 60E15.
\par\noindent\rule{\textwidth}{0.4pt}

\section{Introduction}\label{sec1}

The fundamental concept of entropy, serving as a quantitative measure of disorder or uncertainty, has a rich history originating in thermodynamics. The notion of entropy was introduced by \citet{clausius1850ueber} in relation to the second law of thermodynamics. Later, \citet{boltzmann1872weitere} provided a statistical definition of entropy by linking it to statistical mechanics. A pivotal moment in its evolution was the seminal work of \citet{shannon1948mathematical}, which introduced entropy as a core concept in information theory, precisely quantifying the uncertainty inherent in probability distributions. If $X$ is a nonnegative random variable (RV) with an absolutely continuous cumulative distribution function (CDF) $F(x)$ and probability density function (PDF) $f(x)$, Shannon's entropy, $H(X)$, is defined by
\begin{equation}
H(X) = - \int_0^\infty f(x) \log f(x) \, \mathrm{d}x.
\label{eq:ShannonEntropy}
\end{equation}
Shannon's entropy has gained importance in a wide array of disciplines, including communication theory, statistical inference, and reliability engineering, for assessing the overall information content or unpredictability of a system. For some of the notable contributions across these areas, see \citep{saha2024weighted, sathar2021dynamic, kerridge1961inaccuracy}.

The demand for more flexible measures of uncertainty, capable of capturing different characteristics of probability distributions, led to several important generalizations. Among these, \citet{renyi1961measures} introduced a parameterized family of entropies, providing a broader analytical framework. For a nonnegative RV $X$ with PDF $f(x)$, Rényi's entropy of order $\beta$ is defined as
\begin{equation}
r_\beta(X) = \frac{1}{1 - \beta} \log \left( \int_0^\infty f^\beta(x) \, \mathrm{d}x \right), \quad \beta \in \mathcal{B}=(0,1) \cup (1, \infty).
\label{eq:RenyiEntropy}
\end{equation}
It is easy to see that as the parameter $\beta$ approaches 1, its limit converges to Shannon's entropy
\begin{equation*}
\lim_{\beta \to 1} r_\beta(X) = - \int_0^\infty f(x) \log f(x) \, \mathrm{d}x = H(X).
\end{equation*}
Beyond its theoretical elegance, Rényi entropy has found extensive practical applications across diverse fields due to its tunable parameter $\beta$, which allows it to emphasize different aspects of the probability distribution. In image processing, Rényi entropy is widely employed for tasks such as image segmentation and thresholding, enabling robust partitioning of images based on pixel intensity distributions \citep[see][]{sahoo1988survey, pun1980new}. Its crucial role also extends to image fusion \citep{gabarda2005renyi}, edge detection \citep{el2012edges}, and medical image registration \citep{pan2017analysis}. Furthermore, its utility is observed in signal processing for analysing the complexity of time series data \citep{gabarda2005renyi, hughes2009application}, pattern recognition for feature extraction and classification \citep{pal1988renyi}, and various other areas like ecology, quantum information, and even finance.

Building on the concept of generalized entropies, \citet{tsallis1988possible} proposed another significant extension, known as Tsallis entropy of order $\beta$. This non-extensive entropy has proven particularly valuable in the study of complex systems and non-extensive statistical mechanics. For a continuous random variable $X$, Tsallis entropy is defined as
\begin{equation}
t_\beta(X) = \frac{1}{\beta - 1} \left[ 1 - \int_0^\infty \left( f(x) \right)^\beta \, \mathrm{d}x \right], \quad \beta \in \mathcal{B}.
\label{eq:TsallisEntropy}
\end{equation}
Similar to Rényi's entropy, Tsallis entropy also recovers Shannon's entropy in the limit as $\beta \to 1$. These generalized entropies provide increased versatility in characterizing diverse uncertainty profiles, finding applications in fields ranging from statistical physics and thermodynamics to signal processing. A known relationship between Rényi and Tsallis entropies is given by
\begin{equation}
r_\beta(X) = \frac{1}{1-\beta} \log \left[ 1 + (1-\beta) t_\beta(X) \right].
\label{eq:RelationRenyiTsallis}
\end{equation}
Entropies given in \eqref{eq:ShannonEntropy}, \eqref{eq:RenyiEntropy}, and \eqref{eq:TsallisEntropy} can take negative values for some random variables. Also, approximating them using empirical distributions is challenging. To overcome these limitations, further extending the analysis of uncertainty from a cumulative perspective, \citet{rao2004cumulative} introduced the cumulative residual entropy (CRE). Unlike the density-based Shannon entropy, CRE possesses distinct advantages: it is always nonnegative and can be defined for both discrete and continuous random variables, thereby overcoming cases where Shannon entropy may become negative. For a nonnegative random variable $X$, CRE, denoted by $\xi(X)$, is defined as
\begin{equation*}
\xi(X) = - \int_0^\infty \bar{F}(x) \log \bar{F}(x) \, \mathrm{d}x.
\label{eq:CRE}
\end{equation*}
This measure quantifies the uncertainty contained in the survival function $\bar{F}(x)$ (or equivalently in the cumulative distribution function $F(x)$), offering insights complementary to those provided by the probability density function–based measures. Several extensions and generalizations of CRE have been proposed in the literature, including the generalized cumulative residual entropy introduced by \citet{psarrakos2017generalized}. For more details about CRE, we refer the reader to \citet{toomaj2017some}.

Building upon these foundational concepts, researchers have further extended generalized entropies to the cumulative residual context. A cumulative residual version of Rényi's entropy, denoted as $I_\beta(X)$, was independently introduced by \citet{sunoj2012dynamic} and \citet{abbasnejad2010dynamic}, defined as
\begin{equation}
I_\beta(X) = \frac{1}{1 - \beta} \log \left( \int_0^\infty \bar{F}^\beta(x) \, \mathrm{d}x \right), \quad \beta \in \mathcal{B}.
\label{eq:CRRenyiEntropy_SunojAbbasy}
\end{equation}
Additionally, extensions of Tsallis entropy to the cumulative residual domain have been proposed. \citet{sati2015some} defined a version, $T_\beta^*(X)$, as
\begin{equation*}
T_\beta^*(X) = \frac{1}{\beta-1} \left( 1 - \int_0^\infty \bar{F}^\beta(x) \, \mathrm{d}x \right), \quad \beta \in \mathcal{B}.
\label{eq:CRTsallisEntropy_Sati_Gupta}
\end{equation*}
This was further extended by \citet{rajesh2019some}, who introduced another form of cumulative residual Tsallis entropy, $T_\beta(X)$, given by
\begin{equation}
T_\beta(X) = \frac{1}{\beta - 1} \int_0^\infty \left( \bar{F}(x) - \bar{F}^\beta(x) \right) \mathrm{d}x, \quad \beta \in \mathcal{B}.
\label{eq:CRTsallisEntropy_Sunoj}
\end{equation}
Notice that as $\beta \to 1$, $T_\beta(X) \to \xi(X)$. Extensive research has been conducted on the various properties of these cumulative residual Tsallis entropies, with notable contributions by \citet{toomaj2022some}, \citet{cali2017some}, \citet{mohamed2022cumulative}, \citet{khammar2018weighted}, \citet{sunoj2018quantile} and \citet{cali2017some}.

Despite these developments, the cumulative residual Rényi entropy $I_\beta(X)$ as defined in \eqref{eq:CRRenyiEntropy_SunojAbbasy} exhibits certain limitations. Firstly, it does not consistently converge to the CRE $\xi(X)$ as $\beta \to 1$. It can be easily seen as
\begin{align*}
\lim_{\beta \to 1} I_\beta(X) &= \lim_{\beta \to 1} \frac{1}{1 - \beta} \log \left( \int_0^\infty \bar{F}^\beta(x) \, \mathrm{d}x \right) \\
&= \frac{\int_0^\infty \bar{F}(x) \log \bar{F}(x) \, \mathrm{d}x}{- \int_0^\infty \bar{F}(x) \, \mathrm{d}x} \\
&= \frac{\xi(X)}{\int_0^\infty \bar{F}(x) \, \mathrm{d}x}.
\end{align*}
Here, the denominator $\int_0^\infty \bar{F}(x) \, \mathrm{d}x$ represents the expectation of $X$, $\mathrm{E}(X)$. This expectation can be any nonnegative number and is not necessarily equal to 1. Therefore, $\lim_{\beta \to 1} I_\beta(X) \neq \xi(X)$, unless $\mathrm{E}(X)=1$. Secondly, as stated in Remark 2.1 of \citet{sunoj2012dynamic}, if we consider that $X$ follows an exponential distribution with mean $\frac{1}{2}$, then $I_\beta(X)= \frac{1}{\beta-1}\log(2\beta)$. Clearly, for $0.5 < \beta < 1$, $I_\beta(X)$ is negative. So $I_\beta(X)$ is not always guaranteed to be positive, which can be an undesirable property for an uncertainty measure.

These drawbacks, particularly the non-convergence to CRE and the lack of guaranteed positivity for $I_\beta(X)$, coupled with the established relationship between the original Rényi and Tsallis entropies given by \eqref{eq:RelationRenyiTsallis}, motivate us to define an alternative form of cumulative residual Rényi entropy in this paper. This new information measure aims to overcome the identified limitations while maintaining consistency with established information-theoretic principles.

This paper is organised as follows. Section \ref{Sec:CRRE} introduces the cumulative residual Rényi entropy (CRRE), a new information measure that addresses limitations of existing formulations, where we rigorously define CRRE, establish its mathematical properties, and present illustrative examples. Comparisons with related measures, specifically Shannon entropy, cumulative residual Tsallis entropy (CRTE), and cumulative residual entropy (CRE), are provided, alongside a crucial convergence result demonstrating that CRRE inherits the uniform convergence of survival functions. We also obtain various ordering results. Section \ref{Sec:DCRRE} extends this foundation to the dynamic version (DCRRE). We present its representation via the mean residual life (MRL) function, exploring its properties, and establishing that DCRRE uniquely characterises the underlying distributions. In Section \ref{Sec:EmpiricalEstimation}, a non-parametric estimator of the proposed CRRE is developed, its almost sure convergence is established, and it is validated through Monte Carlo simulation. The validity of CRRE as a tool for quantifying uncertainty is confirmed through its application in chaotic theory and financial data in Section \ref{Sec:Validation}. Finally, Section \ref{Sec:Conclusion} provides the conclusion of the article, summarising the key contributions.

\section{Cumulative Residual R\'enyi Entropy}\label{Sec:CRRE}

In this section, we define the cumulative residual R\'enyi entropy (CRRE). Throughout this paper, we assume that all expectations that appear are finite.

\begin{definition}
For a nonnegative random variable $X$ with an absolutely continuous survival function $\bar{F}(x)$, the cumulative residual R\'enyi entropy of order $\beta \in \mathcal{B}$ is defined as
\begin{align}
R_\beta(X) &=
\begin{cases}
\displaystyle
\frac{1}{1 - \beta} \log \left( 1 - \int_0^\infty \left( \bar{F}(x) - \bar{F}^\beta(x) \right) \, \mathrm{d}x \right), & \quad \beta \in (0,1), \\[3ex]
\displaystyle
\frac{1}{\beta - 1} \log \left( 1 - \int_0^\infty \left( \bar{F}^\beta(x) - \bar{F}(x) \right) \, \mathrm{d}x \right), & \quad \beta \in (1,\infty).
\end{cases}
\label{eq:CRRE}
\end{align}
\end{definition}

It is easy to see that the proposed CRRE is nonnegative for all \( \beta \in \mathcal{B} \) and removes limitations of the existing measure (i.e., $\lim_{\beta \to 1} R_\beta(X) = \xi(X)$).

Table \ref{tab:crre-examples} provides some examples of CRRE for different distributions.

\begin{table}[h]
\centering
\caption{Examples of cumulative residual Rényi entropy (CRRE) for different distributions, for $\beta \in (1,\infty)$.}
\label{tab:crre-examples}
\begin{tabular}{@{}lll@{}}
\toprule
\textbf{Distribution} & \textbf{ $\bar{F}(x)$} & \textbf{$R_\beta(X)$} \\
\midrule
Uniform$(0,a)$ & 
$\bar{F}(x) = 1 - \frac{x}{a}, \quad 0 \leq x \leq a$ & 
$\frac{1}{\beta-1}\log\!\left(1-\frac{a}{\beta+1}+\frac{a}{2}\right)$ \\
[1.5ex]
Exponential$(\lambda)$ & 
$\bar{F}(x) = \mathrm{e}^{-\lambda x}, \quad x \geq 0,\ \lambda>0$ & 
$\frac{1}{\beta-1}\log\!\left(1+\frac{1}{\lambda}-\frac{1}{\lambda\beta}\right)$ \\
[1.5ex]
Pareto$(k,\alpha)$ & 
$\bar{F}(x) = \left(\frac{k}{x}\right)^\alpha, \quad x \geq k$ & 
$\frac{1}{\beta-1}\log\!\left(1+\frac{k}{\alpha-1}-\frac{k}{\alpha\beta-1}\right)$ \\
\bottomrule
\end{tabular}
\end{table}

The proposed measure has the following properties. The proofs of these properties are straightforward so omitted.

\begin{proposition}
\( R_{\beta}(X) = 0 \) if and only if \( X \) is degenerate.
\end{proposition}

\begin{proposition}
Let $X$ be a nonnegative random variable with CDF $F$ and let $\phi$ be a strictly increasing bijective function defined on the support of X. Then the CRRE of \( Y = \phi(X) \) is given by 
\[
R_{\beta}(Y) = \begin{cases}
\frac{1}{1-\beta} \log \left( 1 - \int_{\max(0,  \phi^{-1}(0))}^{\infty} (\bar{F}(t) - \bar{F}^{\beta}(t)) \phi'(t) \mathrm{d}t \right), \quad \forall \beta \in (0,1), \\
\frac{1}{\beta-1} \log \left( 1 - \int_{\max(0,  \phi^{-1}(0))}^{\infty} (\bar{F}^{\beta}(t) - \bar{F}(t)) \phi'(t) \mathrm{d}t \right), \quad \forall \beta \in (1,\infty).
\end{cases}
\]
\end{proposition}

\begin{corollary}
If \( Y = aX + b \) with \( a > 0, b \in \mathbb{R} \), then
\[
R_{\beta}(Y) = \begin{cases}
\frac{1}{1 - \beta} \log \left( 1 - a ( 1 - \mathrm{e}^{(1-\beta) R_{\beta}(X)}) \right), \quad \beta \in (0,1), \\
\frac{1}{\beta-1} \log \left( 1 - a ( \mathrm{e}^{(\beta-1) R_\beta(X) }-1)\right), \quad \beta \in (1,\infty).
\end{cases}
\]
\end{corollary}

This shows that the proposed \(R_\beta(X)\) is a shift-independent measure.
Now, we will provide an upper bound for $R_\beta(X)$ when the underlying random variable has support of finite length.

\begin{theorem}
For a nonnegative random variable \( X \) with the support $[0,a]$ with $a\in(0,\infty)$, the CRRE has the following bound.
\[
0 \leq R_\beta(X) \leq \begin{cases}
\frac{\log\left(1 - a\beta^{\frac{1}{1 - \beta}} + a\beta^{\frac{\beta}{1 - \beta}}\right)}{1 - \beta}, \quad  \beta \in (0,1), \\
\frac{\log\left(1 + a\beta^{\frac{1}{\beta-1}} - a\beta^{\frac{\beta}{1 - \beta}}\right)}{\beta-1}, \quad \beta \in (1,\infty).
\end{cases}
\]
\end{theorem}

\begin{proof}
We define \( f: [0,1] \rightarrow \mathbb{R} \), as
\begin{equation}
f(z) \coloneqq \frac{z-z^{\beta}}{\beta-1},
\label{eq:ConcaveFunction1}
\end{equation} 
where $\mathbb{R}$ denotes the set of real numbers and $\beta \in \mathcal{B}$. It is easy to follow that $f$ attains the maximum at $z_0=\left(\frac{1}{\beta}\right)^{\frac{1}{\beta-1}}$ and \(f(z_0)=\beta^{\frac{\beta}{\beta-1}}\). As $f(z)$ is nonnegative in $[0,1]$, we get 
\[
0 \leq f(z) \leq \beta^{\frac{\beta}{1-\beta}}.
\]
Note that $z= \bar{F}(x)$ gives us $f(z)=\frac{ \bar{F}(x) - \bar{F}^\beta(x) }{\beta-1}$. Using this in the inequality above, we get 
\[
0 \leq \frac{ \bar{F}(x) - \bar{F}^\beta(x) }{\beta-1} \leq \beta^{\frac{\beta}{1-\beta}}.
\]
Further integrating from $0$ to $a$ we obtain
\[
0 \leq \frac{ \int_0^a(\bar{F}(x) - \bar{F}^\beta(x))\, \mathrm{d}x }{\beta-1} \leq a\beta^{\frac{\beta}{1-\beta}}.
\]
Let $\beta \in (0,1)$. Then we have
\[
0 \leq \frac{ -\int_0^a(\bar{F}(x) - \bar{F}^\beta(x) )\, \mathrm{d}x}{1-\beta} \leq a\beta^{\frac{\beta}{1-\beta}}
\]
\[
1 \leq 1-\int_0^a(\bar{F}(x) - \bar{F}^\beta(x)) \mathrm{d}x \leq 1 + a(1-\beta) \beta^{\frac{\beta}{1 - \beta}}
\]
\[
0 \leq \frac{\log\left(1-\int_0^a(\bar{F}(x) - \bar{F}^\beta(x)) \mathrm{d}x\right)}{1-\beta} \leq \frac{\log\left(1 + a(1-\beta) \beta^{\frac{\beta}{1 - \beta}}\right)}{1-\beta}
\]
\[
0 \leq R_\beta(X) \leq \frac{\log\left(1 - a\beta^{\frac{ 1}{1 - \beta}} + a\beta^{\frac{\beta}{1 - \beta}}\right)}{1 - \beta}.
\]
Similarly, for \(\beta \in (1, \infty)\), we get the other bound, completing the result.
\end{proof}

In the next theorem, we will compare the proposed CRRE with the cumulative residual Tsallis entropy (CRTE) $T_\beta(X)$ defined in \eqref{eq:CRTsallisEntropy_Sunoj}.

\begin{theorem}\label{thm:CompareCRRE_CRTE}
For a nonnegative random variable \( X \) with survival function \( \bar{F} \), the inequality
\[
R_\beta(X) < T_\beta(X)
\]
holds for all \( \beta \in \mathcal{B} \), where $R_\beta(X)$ and $T_\beta(X)$ are defined in \eqref{eq:CRRE} and \eqref{eq:CRTsallisEntropy_Sunoj} respectively.
\end{theorem}

\begin{proof}
Let \( \beta \in \mathcal{B} \). Using the $T_\beta(X)$ in \eqref{eq:CRRE} the expression of \( R_\beta(X) \) can be written as
\[
R_\beta(X) = \begin{cases}
\frac{1}{1-\beta} \log\left( 1 + (1-\beta) T_\beta(X) \right), \quad \beta  \in (0,1), \\
\frac{1}{\beta-1} \log\left( 1 + (\beta-1) T_\beta(X) \right), \quad \beta \in (1,\infty).
\end{cases}
\]
For each $\beta \in (0,1)$, $z>0$, define a function \( f: (0,\infty)\rightarrow \mathbb{R} \), as \(f(z) \coloneqq \frac{\log(1+(1-\beta)z)}{1-\beta}-z\). It can be seen that the derivative \(f'(z)=\frac{-(1-\beta)z}{1+(1-\beta)z}\) is negative for all $z \in (0,\infty)$. It follows that $f(z)$ is a strictly decreasing function for all $z \in (0,\infty)$. Further, $f(0)=\lim_{z \to 0} f(z) =0$. Therefore, $f(z)<f(0)=0$.
Taking $z=T_\beta(X)$ in $f(z)$ and using $f(z)<0$ we get,
\[\frac{\log(1 + (1-\beta)T_\beta(X))}{1-\beta} < T_\beta(X) .
\]
So, \(R_\beta(X) < T_\beta(X)\). Similarly, for $\beta \in (1,\infty)$, we can prove \(R_\beta(X) < T_\beta(X)\) by considering a function 
\(
f(z)=\frac{\log(1+(\beta-1)z)}{\beta-1}-z
\) and hence the result follows.
\end{proof}

Now, we will see a relation between $R_\beta(X), T_\beta(X), \text{and } \xi(X)$.

\begin{corollary}\label{cor:R_T_xi}
For any nonnegative random variable \( X \), the following inequality holds:
\begin{equation}
 R_\beta(X) < T_\beta(X) < \xi(X), \quad \forall \beta \in (1,\infty).
\label{eq:CRRE_CRTE_CRE}
\end{equation}
\end{corollary}

\begin{proof}
The result follows from two established inequalities. From Theorem \ref{thm:CompareCRRE_CRTE}, we immediately obtain \( R_\beta(X) < T_\beta(X) \) for all \( \beta > 1 \). Applying Theorem 2.1 of \citet{toomaj2022some}, we get \( T_\beta(X) < \xi(X) \) for all $\beta>1$. Combining these two results through transitivity yields the complete inequality chain \( R_\beta(X) < T_\beta(X) < \xi(X) \).
\end{proof}

For $\beta \in (0,1)$, the inequality in \eqref{eq:CRRE_CRTE_CRE} does not hold. Consider a random variable $X$ that follows \(U(0,a), \hspace{0.1cm} a>0 \) distribution. Then, we have,
\[
\xi(X) = \frac{a}{4}, \quad R_{\beta}(X) = \frac{\log \left(1 - \frac{a}{2} + \frac{a}{\beta+1} \right)}{1-\beta}, \quad
T_{\beta}(X) = \frac{a}{2(\beta+1)}.
\]

\begin{figure}[h]
\centering
\begin{minipage}{0.48\textwidth}
\centering
\includegraphics[width=\linewidth]{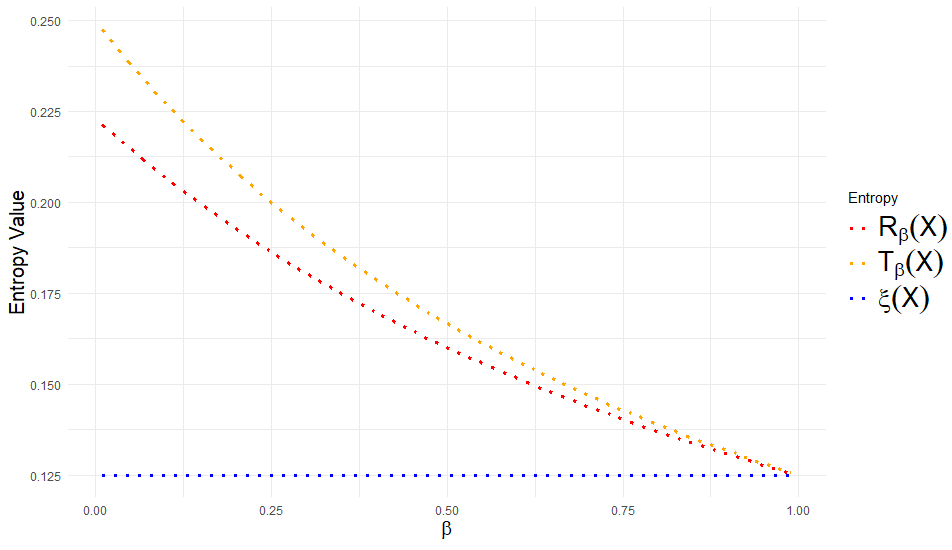}
\caption*{(a) $a = 0.5$}
\end{minipage}
\hfill
\begin{minipage}{0.48\textwidth}
\centering
\includegraphics[width=\linewidth]{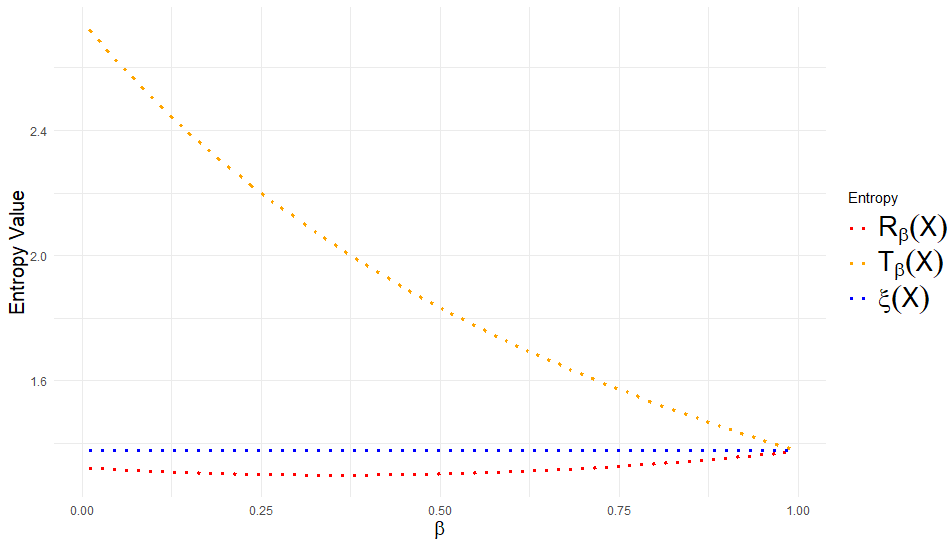}
\caption*{(b) $a = 5.5$}
\end{minipage}
\caption{Comparison of \( \xi(X), R_{\beta}(X), T_{\beta}(X) \) for different values of \( a \).}
\label{fig:ComparisonEntropy}
\end{figure}

We have plotted entropy measures $\xi(X)$, $R_\beta(X)$ and $T_\beta(X)$ for $a=0.5$ and $a=5.5$ in Figure \ref{fig:ComparisonEntropy}. From these, we observe that for \( a = 0.5 \): \( \xi(X) < R_{\beta}(X) < T_{\beta}(X) \) and for \( a = 5.5 \): \( R_{\beta}(X) < \xi(X) < T_{\beta}(X) \). It is easy to see that there is no strict inequality between \( \xi(X) \), \( R_{\beta}(X) \), and \( T_{\beta}(X) \) for \( \beta \in (0,1) \); the ordering depends on the value of \( a \).

Now we discuss the relationship between the proposed CRRE and Shannon entropy defined in \eqref{eq:ShannonEntropy}.

\begin{theorem}
The Shannon Entropy $H(X)$ and the $R_\beta(X)$ shares the following inequality.
$$
R_\beta(X)\geq \begin{cases}\frac{H(X)}{1-\beta} +C_1(\beta),\quad \forall \ \beta \in (0,1), \\
\frac{H(X)}{\beta-1} +C_2(\beta),\quad \forall \ \beta \in (1,\infty).
\end{cases}
 $$
where \(C_1(\beta)= \frac{1}{1-\beta}\int_0^1\log(u^{\beta}-u)\mathrm{d}u \) and \(C_2(\beta)= \frac{1}{\beta-1}\int_0^1\log(u-u^{\beta})\mathrm{d}u \).
\end{theorem}

\begin{proof}
Let $\beta \in (0,1)$. Then by log-sum inequality, we have
\begin{align*}
\int_{0}^{\infty} f(x) \log \left( \frac{f(x)}{f(x) - \left( \bar{F}(x) - \bar{F}^\beta(x) \right)} \right) \mathrm{d}x
&\geq
\log \left( \frac{1}{1 - \int_{0}^{\infty} \left( \bar{F}(x) - \bar{F}^\beta(x) \right) \mathrm{d}x} \right) \\
&= -\log \left( 1 - \int_{0}^{\infty} \left( \bar{F}(x) - \bar{F}^\beta(x) \right) \mathrm{d}x \right) \\
&= -(1 - \beta) R_\beta(X).
\end{align*}
After expanding the left hand side as \( -H(X) - \int_{0}^{\infty} f(x)\log(f(x)-(\bar{F}(x) - \bar{F}^\beta(x))) \, \mathrm{d}x \), we obtain
\begin{equation}
H(X) + \int_{0}^{\infty} f(x)\log\left(f(x)-(\bar{F}(x) - \bar{F}^\beta(x))\right) \, \mathrm{d}x \leq (1-\beta) R_\beta(X).
\label{eq:2.1}
\end{equation}
It is easy to see that for $\beta \in (0,1)$, \( \bar{F}^{\beta}(x)-\bar{F}(x)\leq f(x)+\bar{F}^{\beta}(x)-\bar{F}(x)\). Which further gives us 
\begin{equation}
\int_{0}^{\infty} f(x)\log\left(\bar{F}^\beta(x)-\bar{F}(x)\right) \, \mathrm{d}x \leq \int_{0}^{\infty} f(x)\log\left(f(x)-(\bar{F}(x) - \bar{F}^\beta(x))\right) \, \mathrm{d}x.
\label{eq:2.2}
\end{equation}
From \eqref{eq:2.1} and \eqref{eq:2.2}, we get
\[
H(X) + \int_{0}^{\infty} f(x)\log\left((\bar{F}^\beta(x)-\bar{F}(x))\right) \, \mathrm{d}x \leq (1-\beta) R_\beta(X)
\]
\[
H(X) + \int_{0}^{1} \log\left(u^\beta-u\right) \, \mathrm{d}u \leq (1-\beta) R_\beta(X)
\]
\[
R_\beta(X)\geq \frac{H(X)}{1-\beta} + C_1(\beta),\
\text{where } C_1(\beta)= \frac{1}{1-\beta}\int_0^1\log(u^{\beta}-u)\mathrm{d}u. \]
Now, if we follow the similar steps for $\beta \in (1,\infty)$, then we will get the other bound.
\end{proof}

The forthcoming theorem demonstrates that for independent random variables, the CRRE computed for their sum is greater than the largest CRRE value among the constituent variables.

\begin{theorem}
Let $X$ and $Y$ be independent, nonnegative random variables. Then, for any $\beta \in \mathcal{B}$, we have
\[
R_\beta(X + Y) \geq \max \left\{ R_\beta(X), R_\beta(Y) \right\}.
\]
\end{theorem}

\begin{proof}
Let \( X \) and \( Y \) be independent with respective densities \( f_X, f_Y \) and survival functions \( \bar{F}_X, \bar{F}_Y \). Using the convolution formula, we have
\begin{equation}
\mathbb{P}(X + Y > t) = \int_{-\infty}^{\infty} \bar{F}_X(t - y) f_Y(y) \, \mathrm{d}y = \mathrm{E}_Y \left( \bar{F}_X(t - Y) \right), \quad t \in \mathbb{R}.
\label{eq:surv}
\end{equation}
Consider the function $g(z)=\frac{z-z^{\beta}}{\beta-1}$ defined in \eqref{eq:ConcaveFunction1}. It is concave on $[0,1]$ for all $\beta \in \mathcal{B}$. Applying Jensen’s inequality to \eqref{eq:surv} for the concave function \( g \), we get
\begin{equation}
    g\left( \mathrm{E}_Y \left( \bar{F}_X(t - Y) \right) \right) \geq \int_0^\infty g \left( \bar{F}_X(x) \right) f_Y(t-x) \mathrm{d}x.
    \label{eq:jensen-tbeta}
\end{equation}
Let $\beta \in (0,1)$. Then the CRRE will become
\begin{equation}
 R_\beta(X + Y) = \frac{1}{1-\beta}\log \left( 1+ (1-\beta) \int_0^\infty g \left( \mathrm{E}_Y \left( \bar{F}_X(t - Y) \right) \right)\mathrm{d}t \right) .
\label{eq:TsallisFunctional}
\end{equation}
Using the inequality given in \eqref{eq:jensen-tbeta} and changing the order of integration (using Fubini’s theorem) in \eqref{eq:TsallisFunctional}, we get
\begin{align*}
R_\beta(X + Y)
&= \frac{1}{1-\beta}\log \left( 1+ (1-\beta) \int_0^\infty g \left( \mathrm{E}_Y \left( \bar{F}_X(t - Y) \right) \right)\mathrm{d}t \right) \\
&\geq \frac{1}{1-\beta}\log \left( 1+ (1-\beta) \int_0^\infty \mathrm{E}_Y \left( g \left( \bar{F}_X(t - Y) \right) \right)\mathrm{d}t \right) \\
&= \frac{1}{1-\beta}\log \left( 1+ (1-\beta) \int_{0}^\infty\int_0^t g \left( \bar{F}_X(x)  \right)f_Y(t-x) \mathrm{d}x\mathrm{d}t\right) \\
&= \frac{1}{1-\beta}\log \left( 1+ (1-\beta)\int_{0}^\infty\int_x^\infty g\left( \bar{F}_X(x) \right)f_Y(t-x) \mathrm{d}t\mathrm{d}x\right)\\
&= \frac{1}{1-\beta}\log \left( 1+ (1-\beta)\int_{0}^\infty g\left( \bar{F}_X(x) \right) \left(\int_{x}^\infty f_Y(t-x) \mathrm{d}t\right)\mathrm{d}x\right)\\
&= \frac{1}{1-\beta}\log \left( 1+ (1-\beta)\int_{0}^\infty g\left( \bar{F}_X(x) \right) \left(\int_{0}^\infty f_Y(y) \mathrm{d}y\right)\mathrm{d}x \right) \\
&= \frac{1}{1-\beta}\log \left( 1+ (1-\beta)\int_{0}^\infty g\left( \bar{F}_X(x) \right) \mathrm{d}x \right) \\
&= R_\beta(X).
\end{align*}
Thus, \( R_\beta(X + Y) \geq R_\beta(X) \). By symmetry, one can similarly show \( R_\beta(X + Y) \geq R_\beta(Y) \), and therefore \(R_\beta(X + Y) \geq \max \left\{ R_\beta(X), R_\beta(Y) \right\}\). Following the similar steps for $\beta \in (1,\infty)$, we get the result.
\end{proof}

The next theorem establishes a uniform convergence result, which plays a central role in our study. In particular, it shows that the cumulative residual Rényi entropy $R_\beta(X)$ inherit the convergence behaviour of the underlying survival functions.

\begin{theorem}\label{thm:UniformConvergenceResult}
Let \(\{\bar{F}_n(x): n \geq 1\}\) be a sequence of survival functions corresponding to a sequence of random variables \(\{X_n: n \geq 1\}\). Suppose \(\bar{F}_n(x) \to \bar{F}(x)\) uniformly, where \(\bar{F}(x)\) is the survival function of a random variable \(X\). Then $R_\beta(X_n) \to R_\beta(X)$ uniformly.
\end{theorem}

\begin{proof}
We begin by proving the convergence of \( \xi(\bar{F}_n(x)) \). Since \( \bar{F}_n(x) \to \bar{F}(x) \) uniformly, it follows that for a given $\delta>0$ there exists $n_0 \in \mathbb{N}$ such that $|\bar{F}_n(x) - \bar{F}(x)|<\delta$ for every, $n\geq n_0$ and $x \in (0,\infty)$. We know that a continuous function on a compactly supported domain is uniformly continuous. So the function \( f(x) = \frac{x - x^\beta}{\beta - 1} \) is uniformly continuous on \( [0, 1] \). So given $\epsilon>0$, there exists a $\delta>0$ such that $|f(x)-f(y)|<\epsilon$ whenever $|x-y|<\delta$. For $x=\bar{F}_n(x)$ and $y=\bar{F}(x)$ we get
\[
\left |\frac{\bar{F}_n(x)-\bar{F}_n^{\beta}(x)}{\beta-1}- \frac{\bar{F}(x)-\bar{F}^{\beta}(x)}{\beta-1}\right|<\epsilon,
\]
which implies \( g(\bar{F}_n(x)) \to g(\bar{F}(x)) \) uniformly as $n \to\infty $.
Using the bounded convergence theorem, we get
\[
\lim_{n \to \infty}\int_0^\infty \ g(\bar{F}_n(x)) \, \mathrm{d}x =\int_0^\infty \lim_{n \to \infty} g(\bar{F}_n(x)) \, \mathrm{d}x = \int_0^\infty g(\bar{F}(x)) \, \mathrm{d}x.
\] 
For the case $\beta \in (0,1)$, we have \( R_\beta(X_n) =\frac{1}{1-\beta}\log \left( 1+ (1-\beta) \int_0^\infty g(\bar{F}_n(x)) \, \mathrm{d}x\right) \). Using uniform convergence and the continuity of the logarithm function where \( 1 + (1 - \beta)\int_0^\infty g(\bar{F}_n(x)) \, \mathrm{d}x > 0 \), we obtain
\[
\lim_{n \to \infty} R_\beta(X_n) = \frac{1}{1 - \beta} \log\left(1 + (1 - \beta) \lim_{n \to \infty} \int_0^\infty g(\bar{F}_n(x))\mathrm{d}x\right) = R_\beta(X).
\] 
Similarly, we can prove the result for $\beta \in (1,\infty)$. Hence the result follows.
\end{proof}

\begin{remark}
Following the proof of Theorem \ref{thm:UniformConvergenceResult}, one can easily prove that $\xi(X_n) \to \xi(X)$ and  $T_\beta(X_n) \to T_\beta(X)$ uniformly as $n \to\infty $.
\end{remark}

In the next theorem we will see the concavity of our proposed CRRE with respect to a weighted arithmetic mean of survival functions. The proof is straightforward so omitted.

\begin{theorem} \label{thm:WeightedAM}
Let \( X_1, X_2, \ldots, X_m \) be random variables with survival functions \( \bar{F}_i(x) \) for \( i = 1, 2, \ldots, m \). Let \( X \) be a random variable with survival function
\(
\bar{F}(x) = \sum_{i=1}^m \alpha_i \bar{F}_i(x),
\)
where \( 0 \leq \alpha_i \leq 1 \) and \( \sum_{i=1}^m \alpha_i = 1 \). Then $R_\beta(X) \geq \sum_{i=1}^m \alpha_i R_\beta(X_i)$.
\end{theorem}

\subsection{Ordering Results}

Ordering of random variables provides a rigorous framework for comparing distributions, which is crucial in analyzing key characteristics such as lifetimes, risks, and performance. It offers a systematic way to study the ``ageing'' behaviour of distributions, helping to determine whether the expected future lifetime increases or decreases with age. Moreover, ordering plays an important role in establishing inequalities and identifying relationships among different statistical measures. These aspects make it a powerful tool with broad applications in reliability theory, risk management, and actuarial science. In what follows, we present some fundamental definitions and then a sequence of theorems that formalize the ordering properties relevant to our study. Proofs are omitted for brevity, as they can be followed directly from standard arguments. For the definitions, we refer to \citet{shaked2007stochastic}.

\begin{enumerate}
\renewcommand{\labelenumi}{(\roman{enumi})}

    \item \textbf{Stochastic Order ($\leq_{st}$):}   
    A random variable \(X\) is said to be smaller than \(Y\) in the stochastic order, denoted by \( X \leq_{st} Y \), if \(\bar{F}_X(t) \leq \bar{F}_Y(t), \quad \forall t > 0,\)
    where \(\bar{F}_X(t) = 1 - F_X(t)\) and \(\bar{F}_Y(t) = 1 - F_Y(t)\).

    \item \textbf{Mean Residual Life (MRL) Function:}   
    For a nonnegative random variable \(X\), its Mean Residual Life (MRL) function at time \(t\) is defined as
    \[
    m_X(t) = \mathrm{E}[X - t \mid X > t] 
    = \frac{\int_t^\infty \bar{F}_X(u) \, \mathrm{d}u}{\bar{F}_X(t)}, 
    \quad \bar{F}_X(t) > 0.
    \]

    \item \textbf{MRL Order ($\leq_{mrl}$):}   
    We say that \( X \leq_{mrl} Y \) if 
    \(
    m_X(t) \leq m_Y(t), \quad \forall t \geq 0.
    \)

    \item \textbf{Increasing/Decreasing Mean Residual Life:}   
    A random variable \( X \) is said to have \emph{increasing mean residual life (IMRL)} if
    \(
    m_X(t) = \frac{\int_t^\infty \bar{F}_X(s) \, \mathrm{d}s}{\bar{F}_X(t)}
    \)
    is increasing in \(t\).    
    Similarly, \(X\) is said to have \emph{decreasing mean residual life (DMRL)} if \(m_X(t)\) is decreasing in \(t\).

    \item \textbf{Hazard Rate Order ($\leq_{hr}$):}
    We say that $X \leq_{hr} Y$ if $\bar{F}_Y(t)/\bar{F}_X(t)$ is increasing in $t$.

    \item \textbf{Dispersive Order ($\leq_{disp}$):}   
    Let \( X \) and \( Y \) be continuous random variables with quantile functions 
    \( F_X^{-1}(p) \) and \( F_Y^{-1}(p) \).    
    We say that \( X \leq_{disp} Y \) if
    \[
    F_X^{-1}(p) - F_X^{-1}(q) \leq F_Y^{-1}(p) - F_Y^{-1}(q), 
    \quad \forall \; 0 < p < q < 1.
    \]

 \item \textbf{Failure Rate Classes:}   
A random variable \(X\) is said to have an increasing failure rate (IFR) if its survival function \(\bar{F}_X(x)\) is log-concave, and a decreasing failure rate (DFR) if \(\bar{F}_X(x)\) is log-convex.

\end{enumerate}

\begin{theorem}\label{thm:OrderingResult}
Let \(X\) and \(Y\) be nonnegative random variables with probability density functions \(f_X(x)\) and \(f_Y(x)\), and survival functions \(\bar{F}_X(x)\) and \(\bar{F}_Y(x)\), respectively. For all \(\beta \in \mathcal{B}\), the following hold:
\begin{enumerate}
\renewcommand{\labelenumi}{(\roman{enumi})}
    \item If \(X \leq_{disp} Y\) then \(R_\beta(X) \leq R_\beta(Y)\).
    \item If \(X \leq_{hr} Y\) and either \(X\) or \(Y\) is DFR, then \(R_\beta(X) \leq R_\beta(Y)\).
    \item If \(X \leq_{st} Y\) and \(X \leq_{mrl} Y\), and either \(X\) or \(Y\) is IMRL, then \(R_\beta(X) \leq R_\beta(Y)\).
    \item If \(X \geq_{st} Y\) and \(X \geq_{mrl} Y\), and either \(X\) or \(Y\) is DMRL, then \(R_\beta(X) \geq R_\beta(Y)\).
\end{enumerate}
\end{theorem}

\begin{proof}
Assume that $X$ and $Y$ are nonnegative, absolutely continuous random variables 
with densities $f_X, f_Y$ and survival functions $\bar{F}_X, \bar{F}_Y$, respectively.  
We use the following characterization of the dispersive order \citep[see Section~3.B of][]{shaked2007stochastic}:
\begin{equation*}
    X \leq_{\mathrm{disp}} Y
    \quad \Longleftrightarrow \quad
    f_X\!\left(F_X^{-1}(1-t)\right) \geq f_Y\!\left(F_Y^{-1}(1-t)\right),
    \quad \forall\, t \in (0,1).
    \label{eq:disp_density}
\end{equation*}
Equivalently,
\begin{equation}
    \frac{1}{f_Y(F_Y^{-1}(1-t))} - \frac{1}{f_X(F_X^{-1}(1-t))} \geq 0,
    \quad t \in (0,1).
    \label{eq:inv_density_diff}
\end{equation}
For any absolutely continuous random variable $Z$ with density $f_Z$ and survival $\bar{F}_Z$, 
let $t = \bar{F}_Z(x) = 1 - F_Z(x)$. Then $x = F_Z^{-1}(1-t)$ and $\mathrm{d}t = -f_Z(x)\,\mathrm{d}x$. 
Hence we have
\begin{equation}
    \int_0^\infty \big( \bar{F}_Z(x) - \bar{F}_Z^\beta(x) \big) \mathrm{d}x
    = \int_{t=1}^{t=0} (t - t^\beta) \frac{-\mathrm{d}t}{f_Z(F_Z^{-1}(1-t))}
    = \int_0^1 \frac{t - t^\beta}{f_Z(F_Z^{-1}(1-t))}\,\mathrm{d}t.
    \label{eq:change_var}
\end{equation}
Define
\[
I_Z(\beta) := \int_0^\infty \big( \bar{F}_Z(x) - \bar{F}_Z^\beta(x) \big) \mathrm{d}x.
\]
Then by \eqref{eq:change_var},
\begin{equation}
    I_Y(\beta) - I_X(\beta)
    = \int_0^1 (t - t^\beta)
      \left(
      \frac{1}{f_Y(F_Y^{-1}(1-t))} - \frac{1}{f_X(F_X^{-1}(1-t))}
      \right) \mathrm{d}t.
    \label{eq:I_diff}
\end{equation}

Let $\beta \in (0,1)$.  
For $t \in (0,1)$, we have $t^\beta \geq t$, so $(t - t^\beta) \leq 0$.  
From \eqref{eq:inv_density_diff}, the term in parentheses in \eqref{eq:I_diff} is nonnegative.  
Hence the integrand in \eqref{eq:I_diff} is nonpositive, giving
\[
I_Y(\beta) - I_X(\beta) \leq 0
\quad \text{which further implies}  \quad
I_X(\beta) \geq I_Y(\beta).
\]
Using the definition of the proposed CRRE for $\beta \in (0,1)$,
\[
R_\beta(Z) = \frac{1}{1 - \beta}
\log\!\left( 1 - I_Z(\beta) \right),
\]
we obtain \(1 - I_X(\beta) \leq 1 - I_Y(\beta), \) which implies  \(R_\beta(X) \leq R_\beta(Y)\),
since $\tfrac{1}{1-\beta} > 0$.

Now for the case $\beta > 1$.  
For $t \in (0,1)$, $t^\beta \leq t$, so $(t - t^\beta) \geq 0$.  
From \eqref{eq:I_diff} and \eqref{eq:inv_density_diff}, this yields
\[
I_Y(\beta) - I_X(\beta) \geq 0
\quad \Longrightarrow \quad
I_X(\beta) \leq I_Y(\beta).
\]
For $\beta > 1$, note that
\[
\int_0^\infty \big( \bar{F}_Z^\beta(x) - \bar{F}_Z(x) \big) \mathrm{d}x
= - I_Z(\beta).
\]
Then, by the definition of $R_\beta(Z)$ for $\beta>1$,
\[
R_\beta(Z)
= \frac{1}{\beta - 1}
\log\!\left( 1 - \int_0^\infty (\bar{F}_Z^\beta(x) - \bar{F}_Z(x))\,\mathrm{d}x \right)
= \frac{1}{\beta - 1}
\log\!\left( 1 + I_Z(\beta) \right).
\]
Since $I_X(\beta) \leq I_Y(\beta)$ and $\frac{1}{\beta - 1} > 0$, it follows that
\[
R_\beta(X) \leq R_\beta(Y).
\]
Combining both cases, we conclude that
\[
R_\beta(X) \leq R_\beta(Y), \quad \forall\, \beta \in \mathcal{B}.
\]
which completes the proof of part \textit{(i)}. For the part \textit{(ii)} we refer to the theorem 3.B.20 of \citet{shaked2007stochastic} which states that if \( X \leq_{\text{hr}} Y \) and either \( X \) or \( Y \) is DFR, then \( X \leq_{\text{disp}} Y \). By part \textit{(i)}, this implies \( R_\beta(X) \leq R_\beta(Y) \). Parts \textit{(iii)} and \textit{(iv)} are direct consequences of Theorem 2.2 in \citet{toomaj2022some}.
\end{proof}

\section{Dynamic Cumulative Residual R\'enyi's Entropy}
\label{Sec:DCRRE}

The duration of an event is a critical focus in various scientific disciplines, including reliability, survival analysis, economics, and business. In the context of survival analysis and life testing, information about the current age of the component under consideration is often available. Consequently, this age must be factored in when quantifying uncertainty. The measure $R_\beta(X)$ is, obviously, unsuitable in such situations and must be modified to take the age into account. Therefore, capturing the influence of the age $t$ of an individual or item on the information content of its residual lifetime is essential for many practical applications. In this section, we introduce the dynamic CRRE (DCRRE) as an extension of the proposed CRRE. If $X$ represents the system's lifetime, the residual lifetime (excess of $X$ beyond age $t$) is defined as $X_t = [X-t \mid X > t]$. The survival function for $X_t$ is $\bar{F}(x)/\bar{F}(t)$ for all $x \geq t$. Based on the CRRE defined in \eqref{eq:CRRE} and the survival function of $X_t$, we proceed to define the DCRRE as
\begin{align}
R_\beta(X; t) &=
\begin{cases}
\displaystyle
\frac{1}{1-\beta} \log \left( 1 - \int_t^\infty \left[ \frac{\bar{F}(x)}{\bar{F}(t)} - \left( \frac{\bar{F}(x)}{\bar{F}(t)} \right)^\beta \right] \mathrm{d}x \right), & \ \beta \in (0,1), \\[3ex]
\displaystyle
\frac{1}{\beta-1} \log \left( 1 - \int_t^\infty \left[ \left( \frac{\bar{F}(x)}{\bar{F}(t)} \right)^\beta - \frac{\bar{F}(x)}{\bar{F}(t)} \right] \mathrm{d}x \right), & \ \beta \in (1,\infty).
\end{cases}
\label{eq:DCRRE}
\end{align}

We will now give another representation of DCRRE using the mean residual life function. We introduce a family of random variables $\{X_{\beta} : \beta > 0\}$ through their survival functions, defined as  
\begin{equation*}
    \overline{F}_{\beta}(x) = \big[\overline{F}(x)\big]^{\beta}, 
    \quad x \geq 0,\; \beta > 0,
\end{equation*}
where $\overline{F}(x)$ denotes the baseline survival function and $\beta$ is referred to as the proportional hazard parameter. The mean residual life functions of $X$ and $X_\beta$ are given by  
\begin{equation}
  m(t) = \int_{t}^{\infty} \frac{\bar{F}(x)}{\bar{F}(t)} \, \mathrm{d}x, \quad \text{and} \quad
m_{\beta}(t) = \int_{t}^{\infty} \left(\frac{\bar{F}(x)}{\bar{F}(t)}\right)^{\beta} \, \mathrm{d}x, \quad 
\label{eq:MRL_Defs}
\end{equation}
respectively. Using \eqref{eq:MRL_Defs} in \eqref{eq:DCRRE}, we have
\begin{align}
R_\beta(X; t)&=
\begin{cases}
\displaystyle
\frac{1}{1-\beta} \log \left( 1 - [m(t) - m_\beta(t)] \right), & \beta \in (0,1), \\[3ex]
\displaystyle
\frac{1}{\beta-1} \log \left( 1 - [m_\beta(t) - m(t)] \right), & \beta \in (1,\infty).
\end{cases}
\label{eq:DCRRE_MRL}
\end{align}

Now, we state the following important result which establishes the monotonicity of the proposed measure with respect to time $t$.

\begin{proposition}
Let $X$ be a nonnegative absolutely continuous random variable. Then, $R_\beta(X; t)$ is increasing (decreasing) in $t$ if and only if
\[
\beta m_\beta(t) \geq (\leq)\, m(t), \quad \text{for all } t > 0.
\]
\end{proposition}
It can be shown that the dynamic version of the proposed measure provides a unique characterization of the underlying distribution, as stated in the following theorem.

\begin{theorem}
Let $X$ be a nonnegative, absolutely continuous random variable with cumulative distribution function (CDF) $F(x)$ and hazard rate function $h(t)$. Then $R_\beta(X; t)$ uniquely determines the underlying distribution.
\end{theorem}

\begin{proof}
Let $F_1(t)$ and $F_2(t)$ be two distribution functions corresponding to random variables $X_1$ and $X_2$, with DCRREs $R_\beta(X_1; t)$ and $R_\beta(X_2; t)$, and mean residual life functions $m_1(t)$ and $m_2(t)$, respectively. 

It is known that two distributions are identical if and only if their mean residual life functions are identical, i.e., 
\[
F_1(x) = F_2(x) \quad \text{if and only if} \quad m_1(t) = m_2(t).
\]
This further implies that 
\[
m_{1,\beta}(t) = m_{2,\beta}(t),
\]
where $m_{1,\beta}(t)$ and $m_{2,\beta}(t)$ denote the mean residual life functions of $X_{1,\beta}$ and $X_{2,\beta}$, respectively. Using the definition given in \eqref{eq:DCRRE_MRL}, we then obtain 
\[
R_\beta(X_1; t) = R_\beta(X_2; t).
\]
Hence, $R_\beta(X; t)$ uniquely determines the distribution.
\end{proof}

The next theorem provides a fundamental characterization result: it shows that the exponential distribution is the only distribution for which the proposed dynamic entropy remains constant over time. We will prove this result for $\beta \in (0,1)$ and the proof follows similarly for $\beta \in (1,\infty)$, so omitted here.

\begin{theorem}
Let $X$ be a non-degenerate random variable. Then
\(R_{\beta}(X,t)\) is a constant if and only if $X$ is exponential,
where
\[
R_{\beta}(X,t) = \frac{1}{1-\beta} \log \left( 1 - m(t) + m_{\beta}(t) \right).
\]
\end{theorem}

\begin{proof}
Assume $X$ follows an exponential distribution with rate $\lambda$, i.e.,
\(
\bar{F}(x) = \mathrm{e}^{-\lambda x}, x>0, \lambda > 0.
\)
Then, the functions $m(t)$ and $m_{\beta}(t)$ simplify to
\[
m(t) = \frac{1}{\lambda}, \quad m_{\beta}(t) = \frac{1}{\lambda \beta}.
\]
Substituting these into the expression for $R_{\beta}(X,t)$, we obtain
\[
R_{\beta}(X,t) = \frac{1}{1-\beta} \log \left( 1 - \frac{1}{\lambda} + \frac{1}{\lambda \beta} \right) = C.
\]
which proves the converse part. Now, assume $R_{\beta}(X,t) = C$, which implies
\[
1 - m(t) + m_{\beta}(t) = C',
\]  
for some constant $C'$. Differentiating both sides with respect to $t$, we obtain
\begin{equation}
 m'(t) = m'_{\beta}(t).
 \label{eq:MRL_1}
\end{equation}
From the definitions of $m(t)$ and $m_{\beta}(t)$, we have the identities:
\begin{equation}
\begin{cases}
1 + m'(t) = m(t) \cdot h(t), \\
1 + m'_{\beta}(t) = m_{\beta}(t) \cdot h_{\beta}(t),
\end{cases}
\label{eq:MRL_relationships}
\end{equation}
where, \( h(t) = \frac{f(t)}{\overline{F}(t)} \) is the hazard rate function and \( h_{\beta}(t) = \beta h(t) \).

From \eqref{eq:MRL_relationships}, we deduce
\begin{equation*}
m(t) = \beta \cdot m_{\beta}(t).
\label{eq:m_beta_relation}
\end{equation*}
Differentiating this equation with respect to $t$, we get  
\begin{equation}
 m'(t) = \beta \cdot m'_{\beta}(t).
\label{eq:MRL_2}
\end{equation}
Now, from \eqref{eq:MRL_1} and \eqref{eq:MRL_2}, we conclude that
\[
\beta = 1 \quad \text{or} \quad m'(t) = m'_{\beta}(t) = 0.
\]
Since $\beta \neq 1$, we must have $m'(t) = 0$, implying that $m(t)$ is constant. This is a well-known characterisation of the exponential distribution, which completes the proof.
\end{proof}

\section{Non-Parametric Estimation of \texorpdfstring{\( R_\beta(X) \)}{R\_beta(X)}}
\label{Sec:EmpiricalEstimation}

In this section, we estimate the cumulative residual Rényi entropy $R_\beta(X)$ using its empirical version. Let $X_1, X_2, \ldots, X_n$ be independently and identically distributed (i.i.d.) absolutely continuous nonnegative random variables, with $X_{(1)} \leq X_{(2)} \leq \cdots \leq X_{(n)}$ as their associated order statistics. The empirical survival function corresponding to $\bar{F}(x)$ is given by
\[
\bar{F}_n(x) = \frac{1}{n}\sum_{i=1}^{n} \textbf{1}_{\{X_i > x\}}, \quad x \geq 0,
\]
where $\textbf{1}_{\{X_i > x\}}$ is the indicator function, taking the value 1 if $X_i > x$ and 0 otherwise. In particular, $\bar{F}_n(x)=1-i/n$ for $X_{(i)}\leq x < X_{(i+1)}$.

We now define the empirical cumulative residual Rényi entropy of order \( \beta \in (0,1) \) as
\[
\widehat{R}_\beta(X) = \frac{1}{1 - \beta} \log \left( 1 - \int_0^\infty \left[ \bar{F}_n(x) - \bar{F}_n^{\,\beta}(x) \right] \, \mathrm{d}x \right).
\]
This expression can be simplified using the order statistics and the sample spacings \( U_{i+1} = X_{(i+1)} - X_{(i)} \), for \( i = 1, 2, \ldots, n-1 \), as
\[
\widehat{R}_\beta(X) = \frac{1}{1 - \beta} \log \left( 1 - \sum_{i=1}^{n-1} U_{i+1} \left[ \left(1 - \frac{i}{n} \right) - \left(1 - \frac{i}{n} \right)^\beta \right] \right).
\]
For any \( \beta \in \mathcal{B} \), we have
\begin{equation}
\widehat{R}_\beta(X) 
=\begin{cases}
\displaystyle
\frac{1}{1 - \beta} \log \left( 1 - \sum_{i=1}^{n-1} U_{i+1} \left[ \left(1 - \frac{i}{n} \right) - \left(1 - \frac{i}{n} \right)^\beta \right] \right), & \quad \beta \in (0,1), \\[3ex]
\displaystyle
\frac{1}{\beta-1} \log \left( 1 - \sum_{i=1}^{n-1} U_{i+1} \left[ \left(1 - \frac{i}{n} \right)^\beta -\left(1 - \frac{i}{n} \right)\right] \right), & \quad \beta \in (1,\infty).
\end{cases}
\label{eq:EmpiricalCRRE}
\end{equation}

\begin{theorem}
The empirical cumulative residual Rényi entropy converges almost surely to the true cumulative residual Rényi entropy. Specifically, for any \( \beta \in \mathcal{B} \), as \( n \rightarrow \infty \), we have
\[
\hat{R}_{\beta}(X_n) \rightarrow R_{\beta}(X) \quad \text{a.s.}
\]
\end{theorem}

\begin{proof}
By the Glivenko–Cantelli theorem, the empirical survival function \( \bar{F}_n(x) \) converges uniformly almost surely to the true survival function \( \bar{F}(x) \), i.e.,
\begin{equation*}\label{eq:glivenko}
\sup_{x \in \mathbb{R}} \left| \bar{F}_n(x) - \bar{F}(x) \right| \rightarrow 0, \quad \text{a.s.}
\end{equation*}
For more details, see \citep{gibbons2014nonparametric}.

Since \( \hat{R}_\beta(X_n) \) is a continuous functional of \( \bar{F}_n(x) \), and the integrand is bounded on compact subsets of the support for \( \beta > 0 \), the result follows by the continuous mapping theorem and the bounded convergence theorem.
\end{proof}

Now, we will validate this convergence with the help of examples for various distributions. We will use the theoretical expressions of $R_\beta(X)$ for these distributions given in Table \ref{tab:crre-examples}.

\begin{figure}[h]
\centering
\begin{subfigure}[t]{0.32\textwidth}
    \centering
    \includegraphics[width=\linewidth]{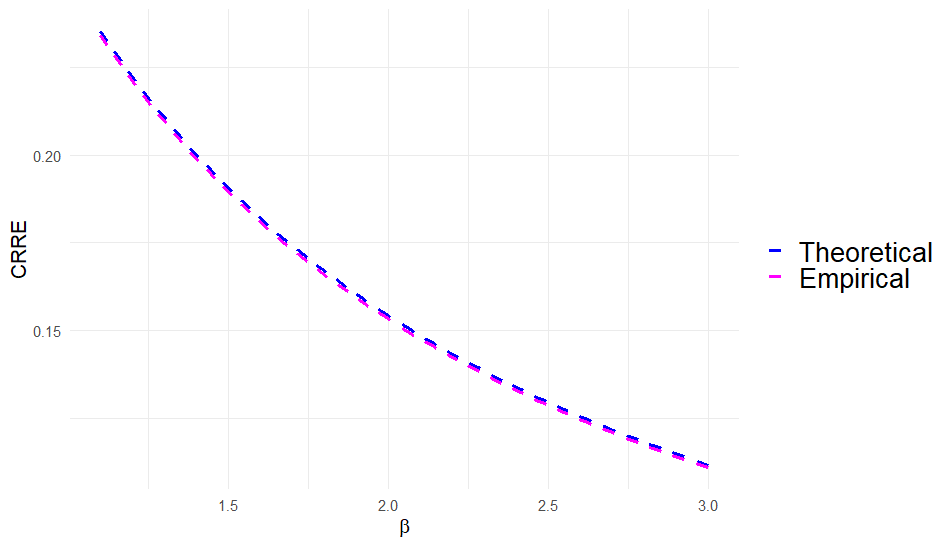}
    \caption{Uniform distribution on \([0,1]\).}
    \label{fig:uniform_crre_beta_gt_1}
\end{subfigure}\hfill
\begin{subfigure}[t]{0.32\textwidth}
    \centering
    \includegraphics[width=\linewidth]{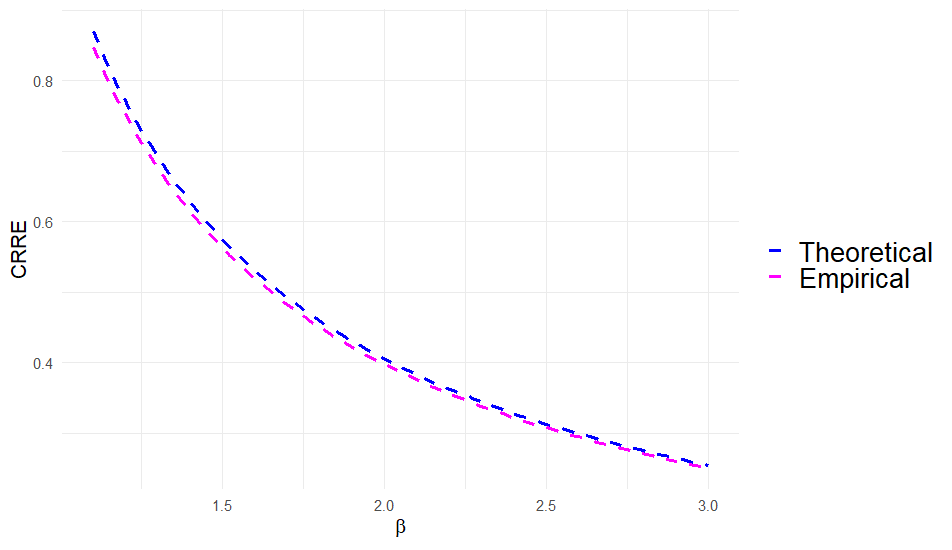}
    \caption{Exponential distribution with mean 1.}
    \label{fig:exponential_crre_beta_gt_1}
\end{subfigure}\hfill
\begin{subfigure}[t]{0.32\textwidth}
    \centering
    \includegraphics[width=\linewidth]{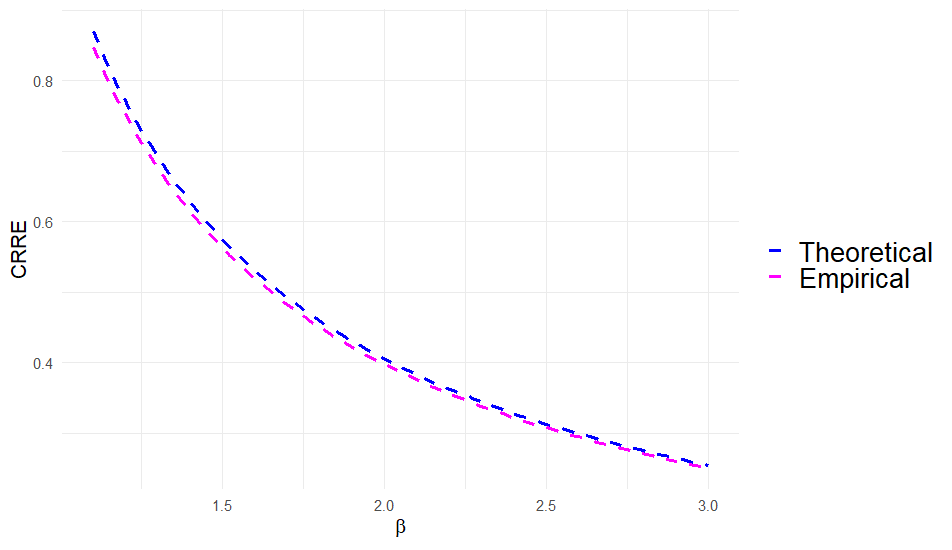}
    \caption{Pareto distribution with $k=1, \alpha=2$.}
    \label{fig:pareto_crre_beta_gt_1}
\end{subfigure}
\caption{Comparison between Theoretical CRRE and Empirical CRRE for various distributions.}
\label{fig:all_crre}
\end{figure}

Figures \ref{fig:uniform_crre_beta_gt_1}, \ref{fig:exponential_crre_beta_gt_1}, and \ref{fig:pareto_crre_beta_gt_1} illustrate the comparison between the theoretical and empirical values of CRRE for the Uniform, Exponential, and Pareto distributions, respectively, for various values of $\beta > 1$.

To validate the convergence of our empirical estimation procedure for the cumulative residual Rényi entropy (CRRE), we generated 10\,000 random samples from each distribution (Uniform on \([0,1]\), Exponential with mean 1, and Pareto with $k=1, \alpha=2$). For each case, the theoretical CRRE was computed using its closed-form expression, while the empirical estimate was obtained numerically. As depicted in Figures \ref{fig:uniform_crre_beta_gt_1}, \ref{fig:exponential_crre_beta_gt_1}, and \ref{fig:pareto_crre_beta_gt_1}, the empirical values for $\beta > 1$ consistently align well with the theoretical ones. This strong agreement confirms the accuracy and consistency of our empirical estimation method across different distributions.

In the next section, we will see the use of empirical CRRE given in \eqref{eq:EmpiricalCRRE} in the validation and application of the proposed CRRE, respectively. For theoretical background crucial to this discussion, we direct the reader to \citep{xiong2019fractional}.

\section{Validation and Application of the Proposed CRRE}
\label{Sec:Validation}

\subsection{Experimental Setup}

All experiments were conducted on a PC equipped with an Intel\textregistered\ Core\texttrademark\ i7-12700 CPU operating at 2.10 GHz and 16 GB of RAM. The system runs on the Windows\textregistered\ 11 Pro operating system (Version 24H2). MATLAB\textregistered\ R2024b was used for the implementation and execution of all algorithms, while RStudio (Version 2024.04.2, Build 764) was employed for statistical computing and visualization tasks. This setup is used for the numerical simulation and data analysis throughout the study.

\subsection{Validity of CRRE with Chaotic Theory}

In this part, we assess the reliability of the proposed measure through simulations on the logistic map \citep{may1976simple}, a classic model for studying chaotic dynamics. The logistic system is defined as $x_{t+1} = a x_t (1 - x_t)$, where $x \in [0,1]$ and $a \in [0,4]$. For our simulations, the initial value is fixed at $x_0 = 0.1$, and the control parameter $a$ is set to $3.2, 3.5, 3.56, 3.8,$ and $4$, which correspond to a period-2, period-4, period-8, and two chaotic series, respectively. This range of parameters allows us to test if the CRRE can effectively distinguish between periodic and chaotic behaviour. 

For each simulation, a total of 20\,000 iterations were performed. The first 1\,000 iterations were discarded to remove transient effects, ensuring that the analysis is based on the system's steady-state behaviour. This results in a fixed sample size of 19\,000 data points for the CRRE calculation, ensuring a fair and consistent comparison across all values of $a$.

\begin{figure}[h]
    \centering
    \includegraphics[width=0.85\textwidth]{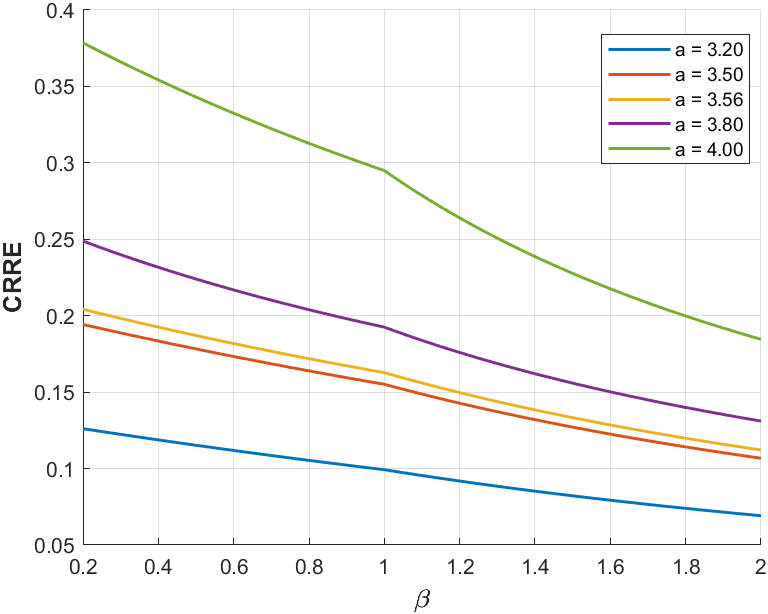}
    \caption{The CRRE of the logistic map under varying parameter values $a$. The sample size is fixed at 19\,000 data points, while the order parameter $\beta$ is varied from $0.2$ to $2$ with increments of $0.02$.}
    \label{fig:logistic_validation}
\end{figure}

Figure \ref{fig:logistic_validation} presents the corresponding CRRE, $R_{\beta}(X)$, for these generated sequences. The results show that CRRE effectively quantifies the uncertainty levels within the logistic map. The CRRE values for the chaotic regimes ($a = 3.8$ and $a = 4.0$) are higher than those for the periodic cases ($a = 3.2, 3.5, 3.56$). This demonstrates that the CRRE correctly identifies and quantifies the increased irregularity and unpredictability of chaotic systems. Furthermore, as the parameter $a$ increases, the CRRE values consistently rise. This is a direct consequence of the logistic map's transition from stable periodic behaviour to increasing complexity and, ultimately, to chaos. The relative ordering of the curves from the lowest for $a = 3.2$ to the highest for $a = 4.0$ directly reflects this increase in system complexity. The green line, representing the fully chaotic case of $a = 4.0$, is the highest curve on the plot for all values of $\beta$. This confirms that the CRRE measure correctly identifies the state of maximum randomness for the logistic map. These observations confirm that CRRE provides a valid and robust quantification of uncertainty, demonstrating its potential for broader applications in complex systems.

\subsection{Application of CRRE in Financial Data}

In this subsection, we demonstrate the applicability of the proposed CRRE measure in analysing financial time series. For this purpose, we consider the daily closing prices of the Dow Jones Industrial Average (DJIA) over the period from October 15, 1997, to October 24, 2024, comprising a total of $T = 6624$ observations, obtained from \href{https://finance.yahoo.com/quote/^DJI/history/}{Yahoo Finance}. The price return series is defined as
\begin{equation}
    r_t = \log(x_t) - \log(x_{t-1}),
\end{equation}
where $x_t$ denotes the closing price at day $t$. Since the return series $\{r_t\}$ may take negative values, a simple transformation is applied to make the sequence nonnegative for entropy evaluation:
\[
    y_t = r_t - \min \{ r_t \}_{t=1}^{T}.
\]
To capture the temporal evolution of uncertainty in the financial market, we compute CRRE using an overlapping sliding window of length $W = 200$ data points with a shift size of $100$ points. This moving-window framework allows us to detect local changes in uncertainty over time rather than relying on a single global estimate.

\begin{figure}[h]
\centering
 \includegraphics[width=0.95\textwidth]{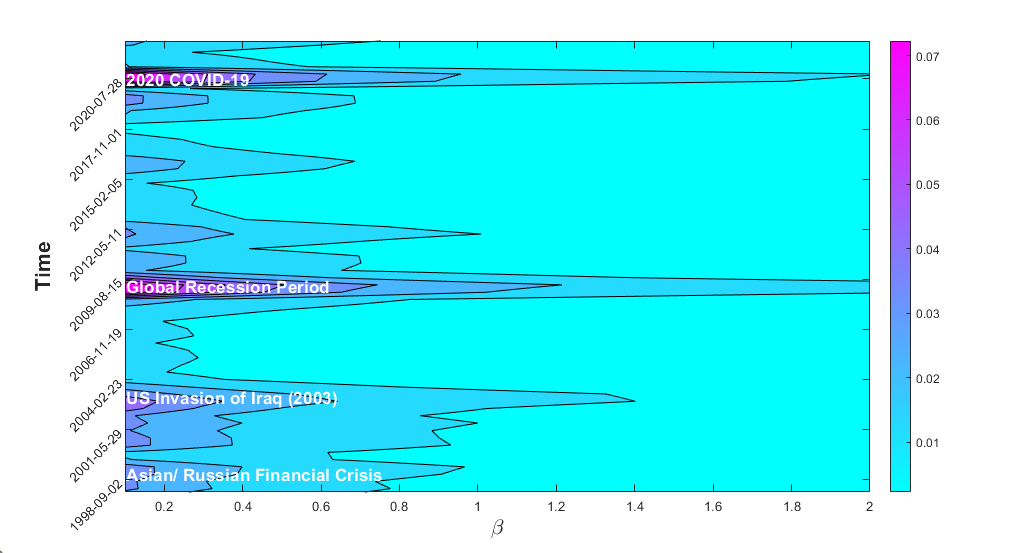}
\caption{Contour representation of CRRE with respect to the order parameter $\beta$ and the time index $t$. The parameter $\beta$ is varied from $0$ to $2$ in steps of $0.02$. The classical CRE corresponds to the special case $\beta = 1$. On the vertical axis, the time index denotes the midpoint of each sliding window.}
\label{fig:finance_application}
\end{figure}

Figure \ref{fig:finance_application} presents a contour plot of CRRE with respect to the order parameter $\beta$ and the time index $t$. The results reveal that CRRE is highly sensitive to variations in financial returns, particularly for $\beta < 0.5$. Several sharp fluctuations are observed that coincide with well-known periods of financial and geopolitical turmoil, including the Asian/Russian financial crisis, the US invasion of Iraq (2003), the 2007–2009 Global Recession, and the COVID-19 pandemic in 2020. This illustrates the ability of CRRE to effectively capture heightened market uncertainty during such events.

Moreover, when $\beta$ is close to one, the variations in CRRE become less visible, and its behaviour is similar to that of the classical CRE. This shows that while CRE gives a basic measure of uncertainty, CRRE, by adjusting the order parameter, can capture finer details and reveal more information in the financial time series. Therefore, CRRE provides a more accurate description of market uncertainty and crisis events than CRE.

\section{Concluding Remarks}
\label{Sec:Conclusion}

In this paper, we addressed the limitations of the existing cumulative residual R\'enyi entropy, $I_\beta(X)$, specifically its issues of non-convergence to the Cumulative Residual Entropy (CRE) and lack of guaranteed positivity. To overcome these drawbacks, we proposed a new information measure cumulative residual R\'enyi entropy. We rigorously defined this new information measure and extensively explored its fundamental mathematical properties, including its intriguing relationship with the existing entropy measures. Various ordering results of the proposed CRRE are obtained. Furthermore, the concept was extended to its dynamic version, DCRRE, and its time-dependent characteristics were thoroughly investigated.
The practical applicability of the proposed information measure was demonstrated by deriving closed-form expressions for various standard probability distributions, such as the uniform, exponential, and Pareto distributions. We also provided a method for its empirical estimation and validated its convergence, demonstrating its utility in real-world scenarios. We applied our measure to chaotic systems and financial data, which confirmed its effectiveness in quantifying uncertainty. This new measure offers a more robust and theoretically consistent tool for quantifying uncertainty in residual lifetimes, providing valuable insights for applications in reliability theory, survival analysis, and broader information-theoretic contexts.

\section*{Disclosure statement}
No potential conflict of interest was reported by the author(s).

\section*{Funding}
Mohd. Arshad gratefully acknowledges the financial support provided by the Anusandhan National Research Foundation (ANRF), India, under the Core Research Grant (Grant No. CRG/2023/001230).

\section*{Data availability statement}
The data utilized in this study were obtained from publicly available sources on Yahoo Finance: \href{https://finance.yahoo.com/quote/^DJI/history/}{Yahoo Finance — DJIA Historical Data}.

\section*{Notes on contributors}
All authors contributed equally to the study.

\bibliography{references}

\end{document}